# ON MINIMAL HYPERSURFACES IN EUCLIDEAN SPACES AND RIEMANNIAN MANIFOLDS

Josef Mikeš[a], Sergey Stepanov[b,c], Irina Tsyganok[c]

[a] *Department of Algebra and Geometry, Palacky University, 17. Listopadu 12, 77146 Olomouc, Czech Republic*

[b] *Department of Mathematics, Russian Institute for Scientific and Technical Information of the Russian Academy of Sciences, 20, Usievicha street, 125190 Moscow, Russia*

[c] *Department of Mathematics, Finance University, 49-55, Leningradsky Prospect, 125468 Moscow, Russia*

**Abstract.** This paper establishes the conditions under which minimal and stable minimal hypersurfaces are characterized as hyperplanes in Euclidean spaces and as totally geodesic submanifolds in Riemannian manifolds.

## 1. Introduction

The study of minimal surfaces has a long and rich history. We recall that an $n$-dimensional minimal submanifold $(M, g)$ in an $m$-dimensional Riemannian manifold $(\bar{M}, \bar{g})$, and particularly in the Euclidean space $\mathbb{R}^m$ for $m > n$, is characterized as a critical point of its volume functional. Moreover, such a submanifold $(M, g)$ is considered stable if the second variation of its volume is always non-negative for any normal deformation with compact support.

Recall that Bernstein showed in 1914 that the graph of any real function $u : \mathbb{R}^3 \to \mathbb{R}$ defined over the Euclidean plane $\mathbb{R}^2$ and which is a minimal surface in $\mathbb{R}^3$ must necessarily be a plane. This theorem was later extended to higher dimensions in various papers (see [10, Chapter 5]). The culminating result asserts that any entire $n$-dimensional minimal graph in $\mathbb{R}^{n+1}$ must be represented by a linear function over $\mathbb{R}^n$, provided that $n \leq 7$. However, for $n \geq 8$, there exist nonlinear entire minimal graphs in $\mathbb{R}^{n+1}$. On the contrary, there exists no compact (without boundary) minimal submanifold in Euclidean space $\mathbb{R}^n$ for any $n \geq 2$. This result, along with other significant findings on minimal submanifolds, can be found in the monograph [1, Chapter 5].

As a natural generalization of the well-known "Bernstein problem," do Carmo and Peng [2] demonstrated that any complete, stable, minimal submanifold $(M, g)$ in $\mathbb{R}^3$ must be planar. Additionally, by employing curvature estimates for minimal hypersurfaces, they further established the following result in [3]: Let $(M, g)$ be an oriented, complete, stable minimal hypersurface in $\mathbb{R}^{n+1}$ such that $\int_M \|A_g\|^2 dv_g < +\infty$, where $A_g$ denotes the shape operator of $(M, g)$ and $\|A_g\|^2 = g(A_g, A_g)$. Then $(M, g)$ is necessarily a hyperplane. This result was later generalized in [5] as follows: Let $(M, g)$ be an $n$-dimensional oriented, complete, stable minimal hypersurface in $\mathbb{R}^{n+1}$ such that $\int_M \|A_g\|^n dv_g < +\infty$. Then $(M, g)$ is also a hyperplane. For further interesting results on the stability of minimal submanifolds, one can refer to the monograph [1, §§5.10, 11.3, 11.4, 16.11], paper [4], and also the current text of the present paper.

We remind the reader that the central question in this theory is to understand the nature of $n$-dimensional minimal surfaces in an $m$-dimensional Riemannian manifold $(\bar{M}, \bar{g})$, especially in Euclidean space $\mathbb{R}^m$ for $m > n$. This question does not have a straightforward answer and, in general, poses significant challenges. In this paper, we aim to provide a partial answer in the specific context of hypersurfaces. We focus on establishing certain results regarding the rigidity of complete, minimal, and stable minimal hypersurfaces. We identify the conditions under which such hypersurfaces are hyperplanes in $\mathbb{R}^{n+1}$ and are totally geodesic submanifolds in an $(n + 1)$-dimensional, $n \geq 2$, Riemannian manifold $(\bar{M}, \bar{g})$.

Obtaining rigidity results for minimal and stable minimal hypersurfaces in Euclidean spaces and Riemannian manifolds is a problem that has attracted much interest in the Geometric Analysis community (see, for example, [5] and [6]). Our methods used in this paper are methods of the Bochner technique (see [7]), which is an important part of the Geometric Analysis. They will allow us to obtain new results in this actively studied topic.

## 2. Complete minimal and stable minimal hypersurfaces in Euclidean space

Recall that a hypersurface in Euclidean space is defined as minimal if it represents a critical point of the area functional. Alternatively, minimal surfaces can be described in terms of their extrinsic curvature. Specifically, the second fundamental form $\varphi$ of a minimal hypersurface $(M, g)$ in $\mathbb{R}^{n+1}$ is a harmonic symmetric bilinear form (see [15, p. 350]) because it is traceless and satisfies the Codazzi equations:

$$(\nabla_X \varphi)(Y, Z) = (\nabla_Y \varphi)(X, Z) \qquad (2.1)$$

for any vector fields $X, Y, Z \in TM$. On the other hand, any $\varphi \in C^\infty(S^2 M)$ satisfying Codazzi equations (2.1) is called the Codazzi tensor (see [18, p. 435]).

In our previous work [14], we proved that on a complete Riemannian manifold $(M, g)$ with non-negative sectional curvature, there does not exist a nonzero harmonic symmetric bilinear form $\varphi \in C^\infty(S^2 M)$ such that $\int_M \|\varphi\|^p dv_g < \infty$ at least one $p \geq 1$. Consequently, there cannot exist a complete minimal hypersurface $(M, g) \subset \mathbb{R}^{n+1}$ with non-negative sectional curvature satisfying $\varphi \in L^p(M, g), p \geq 1$, for its second fundamental form $\varphi$.

Furthermore, it is well known from [16] that the sectional curvature of a hypersurface $(M, g) \subset \mathbb{R}^{n+1}$ is non-negative if and only if the Ricci curvature of $(M, g)$ is also non-negative. In [14], we also proved that a symmetric bilinear form $\varphi \in C^\infty(S^2 M)$ on a Riemannian manifold $(M, g)$ is harmonic if it satisfies the Codazzi equations and has a constant trace. Recall that the mean curvature $H$ of $(M, g)$ is given by $H = \frac{1}{n} \text{trace}_g \varphi$. In particular, the hypersurface $(M, g)$ is minimal if $H \equiv 0$ on $M$ (see [18, p. 39]). Based on the above observations, we formulate the following theorem:

**Theorem 1**. *There is no complete non-compact hypersurface $(M, g) \subset \mathbb{R}^{n+1}$ with constant mean curvature (in particular, no complete minimal hypersurface) and non-negative Ricci curvature satisfying the condition $\int_M \|\varphi\|^p dv_g < \infty$ for the second fundamental form $\varphi$ of the manifold $(M, g)$ and for at least one $p \geq 1$.*

**Remark.** When investigating the general properties of stable minimal hypersurfaces in Euclidean space, the classification of these hypersurfaces has been resolved for dimensions $n = 2,3,4,5$. Specifically, the only orientable complete stable minimal hypersurfaces in $\mathbb{R}^{n+1}$, for $n = 2,3,4,5$, are hyperplanes. For $n \geq 7$, there are examples of orientable complete stable minimal hypersurfaces that are not hyperplanes, and their classification remains an unsolved problem. For instance, in [6], the following theorem was established: A complete, orientable, stable minimal hypersurface in $\mathbb{R}^4$ is a hyperplane.

Let $(M, g)$ be a complete, stable minimal hypersurface in $\mathbb{R}^{n+1}$. Then there exists a positive function $0 < u \in C^\infty(M)$ satisfying the Jacobi equation (see [6])

$$\Delta_g u = -\|A_g\|^2 u, \tag{2.2}$$

where $\Delta_g u := div\,(grad\,u)$ denotes the Laplace–Beltrami operator, $A_g$ represents the shape operator of $(M, g)$ and $\|A_g\|^2 = g(A_g, A_g)$ is the square of the norm of the shape operator. Based on these concepts, we can present the following theorem:

**Theorem 2.** *Let $(M, g)$ be an $n$-dimensional, $n \geq 2$, complete non-compact, stable minimal hypersurface in $\mathbb{R}^{n+1}$ with Ricci curvature $Ric_g$ is bounded below. If there exists a function $0 < u \in L^p(M, g)$ for at least one $0 < p \leq 1$ that solves (2.2), then $(M, g)$ must be a hyperplane. Furthermore, if the volume of $(M, g)$ is infinite, there is no positive function $u \in L^1(M, g)$ that solves equation (2.2).*

*Proof.* From (2.2) we obtain $\Delta_g u \leq 0$. At the same time, a function $u$ on a Riemannian manifold $(M, g)$ is said to be superharmonic if $\Delta_g u \leq 0$. We recall that a manifold $(M, g)$ is called parabolic if any positive superharmonic function on $(M, g)$ is constant, otherwise, it is considered non-parabolic (see [8, p. 164]). If the hypersurface $(M, g)$ is parabolic, it follows from equation (2.2) that it is $\mathbb{R}^n$. It is well-known that $\mathbb{R}^n$ is a parabolic manifold when $n = 2$ and a non-parabolic manifold when $n \geq 3$. Therefore, if $n \geq 3$, then $(M, g) \subset \mathbb{R}^{n+1}$ is non-parabolic. Moreover, there is a generalization of parabolic manifolds known as stochastically complete manifolds. Specifically, any parabolic manifold is stochastically complete, but the converse is not necessarily true (see [8, p. 172]). Recall that if we consider a minimal Wiener process on $(M, g)$, that is, a diffusion process generated by the Laplace–Beltrami operator $\Delta$ with absorption conditions at infinity, then if the probability of absorption at infinity in a finite amount of time is zero, the manifold $(M, g)$ is said to be stochastically complete (see [8], [9], [11]). For instance, $\mathbb{R}^n$ is stochastically complete for all $n \geq 2$ (see [11]). Furthermore, if $(M, g)$ is a stochastically complete manifold, then every positive superharmonic function $u \in L^1(M, g)$ is equal to a constant (see [11]). Consequently, if a hypersurface $(M, g)$ in $\mathbb{R}^{n+1}$ is a stochastically complete Riemannian manifold, it follows from equation (2.2) for each $0 < u \in L^1(M, g)$ we have $A_g = 0$, and thus, consequently, $(M, g)$ is a hyperplane. It is also worth noting that Yau proved in [12] that any complete Riemannian manifold with Ricci curvature bounded below is stochastically complete. Moreover, if the volume of $(M, g)$ is infinite, then under the conditions of Theorem 1, there are no positive superharmonic functions $u \in L^1(M, g)$.

On the other hand, equation (2.2) can be rewritten in the following form:

$$u \, \Delta_g u = - \|A_g\|^2 u^2 \leq 0. \tag{2.3}$$

Recall that Yau established the following result (see [13]): Let $u$ be a smooth function defined on complete Riemannian manifold $(M, g)$ so that $u \geq 0$ and $(p-1) u \, \Delta_g u \geq 0$, where $p$ is a positive number. Then for $p \neq 1$, either $\int_M u^p dv_g = \infty$ or $u$ is a constant. Therefore, for the case $0 < u \in L^p(M, g)$ for at least for one $0 < p < 1$, we deduce from equation (2.3) that $A_g = 0$ and, hence, $(M, g)$ is a hyperplane. The proof is complete.

It is well-known that for any smooth function $u$ on a (not necessarily compact) manifold $(M, g)$ with constant sectional curvature $K$, the Codazzi equations (2.1) imply that $\varphi = \nabla \, du + K \, u \, g$ is a Codazzi tensor (see [18, p. 436]). Moreover, it has been shown that locally, these are the only Codazzi tensors in such manifolds. In particular, from the aforementioned equations, we obtain $trace_g \varphi = \Delta_g u + (nK) u$. In addition, if $trace_g \varphi = 0$, then we have the equation $\Delta_g u + (nK) u = 0$. At the same time, if $\varphi$ is the second fundamental form of a complete, stable minimal hypersurface $(M, g)$ in $\mathbb{R}^{n+1}$, then there exists a positive function $0 < u \in C^\infty(M)$ satisfying the Jacobi equation $\Delta_g u + \|A_g\|^2 u = 0$. From the two preceding equations, we conclude that $K = 1/n \, \|A_g\|^2 \geq 0$. Therefore, we have $\varphi = \nabla \, du$ in the case when $K = 0$, i.e. $(M, g)$ is a Euclidean space form, and $\varphi = K \, g$ in the case when $K > 0$, i.e. $(M, g)$ is a spherical space form.

Suppose $(M, g)$ is a codimension one submanifold in the Euclidean sphere $\mathbb{S}^{n+1}$. The following result is presented in the well-known article [19] on global differential geometry: If $(M, g)$ is a compact minimal hypersurface in the Euclidean sphere $\mathbb{S}^{n+1} \subset \mathbb{R}^{n+2}$ with strictly positive sectional curvature, then it must be an equator of $\mathbb{S}^{n+1}$. We extend this result in the following form:

**Corollary 1.** *Let $(M, g)$ be an n-dimension, $n \geq 2$, compact hypersurface in the Euclidean sphere $\mathbb{S}^{n+1} \subset \mathbb{R}^{n+2}$. If $(M, g)$ has constant mean curvature (in particular, equal to 0) and quasi-positive sectional curvature, then $(M, g)$ is an equator of $\mathbb{S}^{n+1}$.*

*Proof.* To establish this result, we make three key observations. First, if $(M, g)$ is a hypersurface in $\mathbb{S}^{n+1}$ with constant mean curvature, then the second fundamental form $\varphi$ is a Codazzi tensor with constant trace. Second, it was shown in [14] that a symmetric bilinear form $\varphi \in C^\infty(S^2 M)$ on a Riemannian manifold $(M, g)$ is harmonic if it satisfies Codazzi equations (2.1) and has a constant trace. Third, any harmonic symmetric bilinear form $\varphi \in C^\infty(S^2 M)$ on a compact Riemannian manifold $(M, g)$ with quasi-positive sectional curvature must vanish (see [14]). Recall that a manifold with quasi-negative sectional curvature is non-negatively curved and has at least one point where all sectional curvatures are positive. Based on these three observations, we can formulate the above corollary.

**Remark.** The aim of the paper [25] is to prove two results concerning the rigidity of complete, immersed, orientable, stable minimal hypersurfaces: the authors show that they are hyperplane in $\mathbb{R}^4$, while they do not exist in positively curved closed Riemannian $(n + 1)$-manifold when $n \leq 5$; in particular, there are no stable minimal hypersurfaces in $\mathbb{S}^{n+1} \subset \mathbb{R}^{n+2}$ when $n \leq 5$.

## 3. Minimal and stable minimal hypersurfaces in Riemannian manifolds

We recall that one of the most intriguing topics in the calculus of variations within Riemannian geometry is the study of minimal submanifolds $(M, g)$ of a Riemannian manifold $(\bar{M}, \bar{g})$, which arise as critical points of the volume functional. This naturally leads to an investigation of the second variation of the volume of minimal submanifolds, and more specifically, of stable minimal submanifolds (see, for example, [4] and [6]).

To be precise and establish our notation, we note the following. Firstly, we consider smooth, compact, connected hypersurfaces $(M, g) \to (\bar{M}, \bar{g})$ that are isometrically immersed, where $(\bar{M}, \bar{g})$ is a Riemannian manifold of dimension $n + 1 \geq 3$ with metric $\bar{g}$. The induced metric on $M$ is denoted by $g$, and the mean curvature of $(M, g)$ is given by $H = \frac{1}{n} trace A_g$, where $A_g$ denotes the shape operator of $(M, g)$ Moreover, $(M, g)$ is minimal if $H \equiv 0$ (see [18, p. 39]).

It is well-known that the second fundamental form of a hypersurface in a Riemannian manifold of constant sectional curvature is a Codazzi tensor (see [18, p. 436]). Therefore, if $(\bar{M}, \bar{g})$ is a Riemannian manifold of constant sectional curvature, then from Codazzi equations (2.1) we deduce

$$\delta \varphi = -\frac{1}{n} d H. \tag{3.1}$$

At the same time, it is well-known that for any $n$-dimensional $(n \geq 3)$ compact (without boundary) Riemannian manifold $(M, g)$, the algebraic sum Im $\delta^* + C^\infty M \cdot g$ is closed in $S^2 M$, where $\delta^* \theta := \frac{1}{2} L_\xi g$ for the vector field $\xi$ that is dual (by $g$) to the 1-form $\theta$ (see [18, p. 35]). In this case, we have the decomposition

$$S^2 M = (\text{Im } \delta^* + C^\infty M \cdot g) \oplus \left( \delta^{-1}(0) \cap \text{trace}_g^{-1}(0) \right) \tag{3.2}$$

where both factors are infinite dimensional and orthogonal to each other with respect to the $L^2$ inner scalar product $\langle \cdot, \cdot \rangle = \int_M g(\cdot, \cdot) dv_g$ for the canonical measure $dv_g$ of $(M, g)$ (see [18, p. 130]). It's obvious that the second factor $\delta^{-1}(0) \cap \text{trace}_g^{-1}(0)$ of (3.3) is the space of TT-tensors.

**Remark.** Here we recall that a symmetric divergence free and traceless covariant two-tensor is called *TT*-tensor. As a consequence of a result of Bourguignon-Ebin-Marsden (see [18, p. 132]) the space of *TT*-tensors is an infinite-dimensional vector space for any closed Riemannian manifold $(M, g)$. Such tensors are of fundamental importance in stability analysis in General Relativity and in Riemannian geometry (see, for instance, [18, p. 346-347]).

Considering the above, we conclude that the second fundamental form of $\varphi$ has the following $L^2$-orthogonal decomposition (see also formula (3.2))

$$\varphi = \left(\frac{1}{2} L_\xi g + \lambda\, g\right) + \varphi^{TT} \qquad (3.3)$$

for some vector field $\xi \in C^\infty TM$, *TT*-tensor $\varphi^{TT} \in C^\infty(S^2)$ and scalar function $\lambda \in C^\infty M$. If we applying the operator $trace_g$ to both sides of (3.3), we obtain

$$n\, H = -\delta\theta + n\, \lambda \qquad (3.4)$$

where $\theta^\# = \xi$ and $\delta$ is the divergence (see [18, p. 35]). In this case, (3.4) can be rewritten in the form

$$\varphi_0 = S\theta + \varphi^{TT} \qquad (3.5)$$

where $\varphi_0 = \varphi - H\, g$ is the traceless part of the second fundamental form $\varphi$ and $S\theta = \frac{1}{2} L_\xi g + \frac{1}{n} \delta\theta\, g$ is the Cauchy-Ahlfors operator. Next, applying $S^* := \delta$ to both sides of (3.5), we obtain

$$S^* S\, \theta = \delta\varphi_0 \qquad (3.6)$$

for the Ahlfors Laplacian $S^*S$ (see [17]). Combining (3.1) and (3.6), we obtain

$$S^* S\, \theta = -\frac{1}{n} dH. \qquad (3.7)$$

Therefore, if the mean curvature $H$ of $(M, g)$ is constant, then $S^*S\, \theta = 0$ since $\langle S^*S\, \theta, \theta\rangle = \langle S\, \theta, S\theta\rangle \geq 0$ (see [17]). In this case, from (3.6) we obtain $\varphi = H\, g + \varphi^{TT}$. The converse is also true. We are now ready to formulate our result, which generalizes Theorem 5.4.2 from [20].

**Theorem 3.** *Let $(M, g)$ be a compact hypersurface in a Riemannian manifold $(\bar{M}, \bar{g})$ with constant sectional curvature and $\dim \bar{M} \geq 4$. Then the mean curvature $H$ of $(M, g)$ is constant if and only if the second fundamental form $\varphi$ of $(M, g)$ admits the following $L^2$-orthogonal decomposition*

$$\varphi = H\, g + \varphi^{TT} \tag{3.8}$$

*for some TT-tensor $\varphi^{TT}$. In particular, if $(M, g)$ is a compact minimal hypersurface in a Riemannian manifold $(\bar{M}, \bar{g})$ with constant sectional curvature, then $\varphi = \varphi^{TT}$.*

From equation (3.8) and Codazzi equations (2.1), it follows that $\varphi^{TT}$ is a traceless Codazzi tensor, meaning that $\varphi^{TT}$ is a harmonic form. At the same time, we proved in [21] that any such $TT$-tensor must be the zero tensor on a compact Riemannian manifold $(M, g)$ with quasi-positive sectional curvature. In particular, if $(M, g)$ is a hypersurface in a Riemannian manifold $(\bar{M}, \bar{g})$, then it is totally umbilical (see [18, p. 39]) with constant mean curvature. In this case, it follows from (3.8) that $(M, g)$ is a Riemannian manifold of constant sectional curvature. Moreover, this curvature must be positive, given that we assumed that $(M, g)$ to be a manifold of quasi-positive sectional curvature.

A Riemannian manifold with constant sectional curvature is called space form (see [15, p. 20]). Therefore, in our case, $(M, g)$ must be a spherical space form since it has positive constant curvature. In particular, if $(M, g)$ a simply connected manifold, then $M = \mathbb{S}^n$ (see [15, pp. 200-201]). Using this fact, along with the results from Theorem 3, we arrive at the following conclusion:

**Corollary 2**. *Let $(\bar{M}, \bar{g})$ be an $n$-dimensional Riemannian manifold with constant sectional curvature, where $n \geq 4$, and let $(M, g)$ be an isometrically immersed hypersurfaces $(M, g) \to (\bar{M}, \bar{g})$ with constant mean curvature. If $(M, g)$ has quasi-positive sectional curvature, then $(M, g)$ is a spherical space form. In addition, if $(M, g)$ is a simply connected manifold, then it is the Euclidean sphere $\mathbb{S}^n$.*

**Remark**. We present here an alternative proof of Corollary 2. It is well-known the following theorem (see [18, p. 436]): Every Codazzi tensor $\varphi$ with constant trace on a compact Riemannian manifold $(M, g)$ with non-negative sectional curvature is parallel. Moreover, if the sectional curvatures of $(M, g)$ are positive at some point, then $\varphi$ is a constant multiple of $g$. Therefore, if $(M, g)$ is a hypersurface in a Riemannian manifold $(\bar{M}, \bar{g})$ with constant sectional curvature, then it is totally umbilical with constant mean curvature.

Consider $(M, g)$ to be a hypersurface of an $(n + 1)$-dimensional Riemannian manifold $(\bar{M}, \bar{g})$ with positive constant sectional curvature. In such a case, its second fundamental form is a Codazzi tensor (see [18, p. 436]). Moreover, if the mean curvature of $(M, g)$ is constant (and in particular, if it is zero), then the second fundamental form of $(M, g)$ becomes a harmonic bilinear form (see [10]). Using this fact, along with the results from [14], we arrive at the following conclusion:

**Corollary 3.** *Let $(M, g)$ be a complete non-compact hypersurface in an $(n + 1)$−dimensional Riemannian manifold $(\bar{M}, \bar{g})$ with positive constant sectional curvature. If the sectional curvature of $(M, g)$ is non-negative, its mean curvature is constant (in particular, equal to 0) and its second fundamental form $\varphi$ satisfies the condition $\int_M \|\varphi\|^p \, dv_g < \infty$ for at least one $p \geq 1$, then $(M, g)$ must be a spherical space form. In particular, if $(M, g)$ is a simply connected manifold, then it is the Euclidean sphere.*

The following theorem is generalized a Kobayashi's statement (see [20]) for compact minimal hypersurface in a Riemannian manifold of constant curvature.

**Theorem 4**. *Let $(M, g)$ be an n-dimensional complete manifold of finite volume that is minimally immersed in an $(n + 1)$-dimensional Riemannian manifold $(\bar{M}, \bar{g})$ with constant curvature $C$. If $(M, g)$ is not totally geodesic and $\|A_g\|^2 \leq n\,C$, where $A_g$ represents the shape operator of $(M, g)$, then $\|A_g\|^2 = n\,C$ and the second fundamental form $\varphi$ of $(M, g)$ is parallel with respect to the Levi-Civita connection $\nabla$ of $(M, g)$.*

*Proof.* Let $(M, g)$ be an $n$-dimensional complete manifold minimally immersed in an $(n + 1)$-dimensional Riemannian manifold $(\bar{M}, \bar{g})$ of constant curvature $C$. Whether $(M, g)$ is compact or not, according to [20, Eq. (3.10) and Eq. (3.12)], we have the inequality

$$\tfrac{1}{2}\Delta_g \|A_g\|^2 \geq \|\nabla A_g\|^2 + \|A_g\|^2 \left(n\,C - \|A_g\|^2\right) \tag{3.9}$$

for the shape operator $A_g$ of $(M, g)$. Therefore, if we suppose that $(M, g)$ is not totally geodesic and $\|A_g\|^2 \leq n\,C$, then from (3.9) we conclude that $\|A_g\|^2$ is a non-

negative subharmonic function bounded above. Note that the definition of parabolicity allows an equivalent characterization: a complete manifold $(M, g)$ is parabolic if every bounded subharmonic function on $(M, g)$ is a constant (see [8, p. 164], [26, p. 66]). For example, a complete Riemannian manifold $(M, g)$ of finite volume is a parabolic manifold (see [10]). Therefore, if $(M, g)$ has finite volume, then $\|A_g\| = const$. In this case, from (2.5) we obtain $\|A_g\|^2 = n\,C$. The proof is complete.

The famous Eisenhart theorem (see [27, p. 303]) states the following: If $\varphi$ is a parallel symmetric two-tensor on $(M, g)$, then for each point $x \in M$, there is some neighborhood $U \subset M$ where $\varphi = \lambda_1 g_1 + \cdots + \lambda_r g_r$, and locally $(M, g)$ admits a Riemannian direct product structure $(U, g_{|U}) = (U_1, g_1) \times \ldots \times (U_r, g_r)$. Here, the coefficients $\lambda_1, \ldots, \lambda_r$ are constants and $(U_1, g_1), \ldots, (U_r, g_r)$ are Riemannian manifolds of dimensions $n_1 \geq 1, \ldots, n_r \geq 1$, respectively, with $n_1 + \cdots + n_r = n$ for some $r$ (see [27]). In this case the conditions $trace_g \varphi = 0$ and $\|A_g\|^2 = n\,C$ can be rewritten as $n_1 \lambda_1 + \cdots + n_r \lambda_r = 0$ and $n_1 \lambda_1^2 + \cdots + n_r \lambda_r^2 = n\,C$, respectively. Furthermore, according to [20], it is known that if $(M, g)$ is a minimal (but not totally geodesic) hypersurface immersed in an $(n + 1)$-dimensional Riemannian manifold $(\bar{M}, \bar{g})$ of constant curvature 1 satisfying $\|A_g\|^2 = n$, then $r = 2$. In this case, through direct calculations, we deduce the equalities $\lambda_1 = \sqrt{(n-m)/m}$ and $\lambda_2 = -\sqrt{n/(n-m)}$ for $m = \dim U_1 \geq 1$ and $n - m = \dim U_2 \geq 1$. Therefore, we have (see also [20])

**Corollary 4**. *If $(M, g)$ is a minimal (not totally geodesic) submanifold of a Riemannian manifold $(\bar{M}, \bar{g})$ with constant curvature $C = 1$ such that $\|A_g\|^2 = n$, then for each point $x \in M$ there is some neighborhood $U \subset M$ such that $(M, g)$ is locally a Riemannian direct product $(M, g) \supset (U, g_{|U}) = (U_1, g_1) \times (U_2, g_2)$ for Riemannian manifolds $(U_1, g_1)$ and $(U_2, g_2)$ of constant curvatures and dimensions $m \geq 1$ and $n - m \geq 1$, respectively, and $\varphi = \lambda_1 g_1 + \lambda_2 g_2$ with $\lambda_1 = \sqrt{(n-m)/m}$ and $\lambda_2 = -\sqrt{n/(n-m)}$.*

**Remark**. According to [20], the $n$-dimensional minimal submanifold $(M,g)$ in a unit sphere $\mathbb{S}^n$ satisfying the condition $\|A_g\|^2 = n$ has the following local form $\mathbb{S}^m\left(\sqrt{\frac{m}{n}}\right) \times \mathbb{S}^{n-m}\left(\sqrt{\frac{n-m}{n}}\right)$.

Next, we consider complete, connected, isometrically immersed hypersurfaces $(M,g) \to (\bar{M}, \bar{g})$, where $(\bar{M}, \bar{g})$ is a Riemannian manifold of dimension $n+1 \geq 3$. To proceed with further statements, we recall that a minimal hypersurface $(M,g)$ in a Riemannian manifold $(\bar{M}, \bar{g})$ is considered stable if and only if the second variation of its volume functional is non-negative for all compactly supported deformations (see, for example, [1, §5.2], [12]). According to [6], the stability of $(M,g)$, associated with the non-negativity of the second variation, is equivalent to the non-positivity of the Jacobi operator

$$L_g := \Delta_g + \|A_g\|^2 + Ric_{\bar{g}}(v,v),$$

where $Ric_{\bar{g}}(v,v)$ is the Ricci curvature of $(\bar{M}, \bar{g})$ in the direction of the unite normal vector field $v$ to $(M,g)$ at every its point. Therefore, we have

$$\Delta_g u \leq -\left[\|A_g\|^2 + Ric_{\bar{g}}(v,v)\right]u \tag{3.11}$$

for any $u \in C^\infty(M)$. In this case, if $Ric_{\bar{g}}(v,v) \geq 0$ at every point of $M$, then from inequality (3.11) we obtain $\Delta_g u \leq 0$ and, hence, $u$ is a superharmonic function. If a hypersurface $(M,g)$ is a parabolic, complete Riemannian submanifold of $(\bar{M}, \bar{g})$ and there exists a function $0 < u \in L^1(M,g)$ satisfying inequality (3.11), then from (3.11) we deduce that $A_g = 0$ and $Ric_{\bar{g}}(v,v) = 0$ at every point of $M$. In this case, $(M,g)$ is a totally geodesic submanifold of $(\bar{M}, \bar{g})$ This leads to the following straightforward theorem.

**Theorem 5**. *Let $(M,g)$ be a connected, complete, stable minimal hypersurface in a Riemannian manifold $(\bar{M}, \bar{g})$ such that the Ricci curvature of $(\bar{M}, \bar{g})$ in the direction of the unite normal vector field to $(M,g)$ is non-negative at every point of $(M,g)$. If there exists a function $0 < u \in L^1(M,g)$ satisfying the inequality $L_g u \leq 0$ and $(M,g)$ is parabolic, then it must be a totally geodesic submanifold in $(\bar{M}, \bar{g})$.*

**Remark.** This theorem complements the following result (see [4, Theorem 1.1]): Let $(M, g)$ be a complete minimal hypersurface in a manifold $(\bar{M}, \bar{g})$ with nonnegative Ricci curvature $Ric_{\bar{g}}$. If $(M, g)$ is parabolic, then $(M, g)$ must be totally geodesic in $(\bar{M}, \bar{g})$. Moreover, the Ricci curvature $Ric_{\bar{g}}(v, v)$ of $(\bar{M}, \bar{g})$ in the normal direction also vanishes, and $(M, g)$ must has a nonnegative scalar curvature. Note that this can be explained as follows: the non-negativity of the scalar curvature of $(M, g)$ follows from the Gauss curvature equation, combined with the assumption that $(\bar{M}, \bar{g})$ has non-negative Ricci curvature $Ric_{\bar{g}}$.

Furthermore, by the Cheng-Yau theorem, a complete manifold with finite volume is parabolic (see [10]). Consequently, a complete Riemannian manifold of finite volume does not admit non-constant positive superharmonic functions. Moreover, if $0 < u \in C^{\infty}(M)$ and $Ric_{\bar{g}}(v, v) \geq 0$ at every point of $M$, then form (3.8), we conclude that $u$ must be constant, implying $A_g = 0$ and $Ric_{\bar{g}}(v, v) = 0$ at every point of $M$. In this scenario, $(M, g)$ is a totally geodesic submanifold of $(\bar{M}, \bar{g})$. With this in mind, we can formulate the following statement.

**Corollary 5.** *Let $(M, g)$ be a complete, stable minimal hypersurface in a Riemannian manifold $(\bar{M}, \bar{g})$ such that the Ricci curvature of $(\bar{M}, \bar{g})$ in the direction of the unite normal vector field to $(M, g)$ is non-negative at every point of $(M, g)$. If there exists a function $0 < u \in C^{\infty}(M, g)$ satisfying the inequality $L_g u \leq 0$ and the volume of $(M, g)$ is finite, then $(M, g)$ is totally geodesic submanifold in $(\bar{M}, \bar{g})$.*

**Remark.** Insightful details on the stability of totally geodesic submanifolds can be found in the monograph [1, §11.3].

Furthermore, we draw attention here to the elegant Cheng-Yau theorem (see [12, Theorem 1 and Corollary 1]): If, on a complete Riemannian manifold $(M, g)$, the volume $V_R$ of a geodesic ball of radius $R$ with fixed center satisfies the inequality $V_R \leq CR^2$, then every negative subharmonic function on $(M, g)$ must be constant. By definition, a function $u$ is termed subharmonic if $-u$ is superharmonic (see [8, p. 150]). Thus, if $(M, g)$ is a connected, complete, stable minimal hypersurface in $(\bar{M}, \bar{g})$ that meets the condition of the Cheng-Yau theorem and the criteria stated in the previous theorem, then from (3.8), it follows that $A_g = 0$ and, consequently, $(M, g)$ is totally geodesic. In light of the above, we can formulate the following statement.

**Theorem 6**. *Let $(M, g)$ be a connected, complete, stable minimal hypersurface within a Riemannian manifold $(\overline{M}, \bar{g})$ wehere the Ricci curvature of $(\overline{M}, \bar{g})$ in the direction of the unite normal vector field to $(M, g)$ is non-negative at every point of $M$. If the volume $V_R$ of geodesic ball of radius $R$ centered at fixed point satisfies $V_R \leq CR^2$, and there exists a function $0 < u \in C^\infty(M, g)$ that meets the inequality $L_g u \leq 0$, then $(M, g)$ must be a totally geodesic submanifold in $(\overline{M}, \bar{g})$.*

The main result of [11] can be summarized as follows: If on a complete Riemannian manifold $(M, g)$, the volume $V_R$ of a geodesic ball of radius $R$ centered at fixed point satisfies the inequality $V_R \leq e^{CR^2}$, then every nonnegative superharmonic function in $L^1(M, g)$ is a constant. In this scenario, $(M, g)$ is parabolic. Thus, we derive the following statement.

**Theorem 7**. *Let $(M, g)$ be a connected, complete, stable minimal hypersurface in a Riemannian manifold $(\overline{M}, \bar{g})$, where the Ricci curvature of $(\overline{M}, \bar{g})$ in the direction of the unite normal vector field to $(M, g)$ is non-negative at every point of $M$. $L^1(M, g)$ Suppose that the volume $V_R$ of geodesic ball of radius $R$ with a fixed center on $(M, g)$ satisfies the inequality $V_R \leq e^{CR^2}$. If there exists a smooth function $0 < u \in L^1(M, g)$ that satisfies the inequality $L_g u \leq 0$, then $(M, g)$ must be a totally geodesic submanifold in $(\overline{M}, \bar{g})$.*

It is evident that the analogue of Theorem 2 also holds in this case.

**Theorem 8**. *Let $(M, g)$ be a complete, stable minimal hypersurface in a Riemannian manifold $(\overline{M}, \bar{g})$ with the Ricci curvature of $(\overline{M}, \bar{g})$ in the direction of the unite normal vector field to $(M, g)$ being non-negative at every point of $M$. If there exists a smooth function $0 < u \in L^p(M, g)$ that satisfies the inequality $L_g u \leq 0$ for at least one $0 < p \leq 1$, then $(M, g)$ is totally geodesic submanifold in $(\overline{M}, \bar{g})$.*

**Remark.** If $(M, g)$ be an oriented, complete non-compact, stable minimal hypersurface withing a Riemannian manifold $(\overline{M}, \bar{g})$ with non-negative sectional curvatures, then $(M, g)$ has an infinite volume (see [23]). Furthermore, if the volume of $(M, g)$ is infinite, then there cannot be a positive function $u \in L^1(M, g)$ that is satisfies inequality (3.8) with $Ric_{\bar{g}}(v, v) \geq 0$ at every point of $M$.